\documentclass[numreferences]{kluwer}
\usepackage{amssymb}

\newtheorem{prop}{Proposition}
\newtheorem{lem}{Lemma}
\newtheorem{cor}{Corrolary}

\begin{document}
\begin{opening}

\title{\center{Multiple Gamma Function and Its Application to
Computation of Series}}

\author{\centerline{V. S. Adamchik}}


\institute{\center{
Carnegie Mellon University,\\
Pittsburgh, USA}}


\runningtitle{Multiple Gamma Function}
\runningauthor{V. Adamchik}
\begin{abstract}
The multiple gamma function $\Gamma_n$, defined by a
recurrence-functional equation as a generalization of the Euler
gamma function, was originally introduced  by Kinkelin, Glaisher,
and Barnes around 1900. Today, due to the pioneer work of Conrey,
Katz and Sarnak, interest in the multiple gamma function has been
revived. This paper discusses some theoretical aspects of the
$\Gamma_n$ function and their applications to summation of series
and infinite products.
\end{abstract}


\classification{2000 Mathematics Classification}{Primary 33E20,
33F99, 11M35, 11B73}

\end{opening}

\section{Introduction}
The Hurwitz zeta function, one of the fundamental transcendental
functions, is traditionally defined (see \cite{Magnus}) by the
series

\begin{equation}
\label{ZetaDef}
 \zeta (s,z) =\sum _{k=0}^{\infty
}\frac{1}{{{(k+z)}^s}}, \quad \Re(s)>0.
\end{equation}
It admits an analytic continuation to the entire complex plane
except for the simple pole $s=1$. The Riemann zeta function $\zeta
(s)$ is a special case of $\zeta (s,z)$

\[\zeta (s,1) =\zeta (s). \]
\noindent
The Hurwitz function has quite a few series and integral
representations (see \cite{Magnus,Moll02}). The most famous is the
Hermite integral:

\begin{equation}
\label{Hermite}
\begin{array}{c}\displaystyle
 \zeta (s,z)
=\frac{{z^{-s}}}{2}+\frac{{z^{1-s}}}{s-1} + 2\int
_{0}^{\infty}\frac{\sin\big(s\arctan\big(\frac{x}{z}\big)\big)}
{{{({x^2}+{z^2})}^{s/2}}({e^{2 \pi  x}}-1)}\> dx, \\
\noalign{\vskip 1.0pc} \displaystyle
 s\neq 1,\Re(z)>0
\end{array}
\end{equation}
from which one can deduce many fundamental properties of the
Hurwitz function, including the asymptotic expansion at infinity:

\[\zeta
(s,z)=\frac{{z^{1-s}}}{s-1}+\frac{{z^{-s}}}{2}+\sum
_{j=1}^{m-1}\frac{{B_{2 j}}\, \Gamma (2 j+s-1) }{(2 j)!\, \Gamma
(s)}{z^{-2j-s+1}} + O\big(\frac{1}{{z^{ 2m+s+1}}}\big)
\]
The Hurwitz function is closely related to the multiple gamma
function ${\Gamma}_n(z)$ defined as a generalization of the
classical Euler gamma function $\Gamma (z)$, by the following
recurrence-functional equation (for references and a short
historical survey see \cite{AMC03,Vardi88}):

\begin{eqnarray}
\label{GDef} {\Gamma}_{n+1}(z+1) &=&
\frac{{\Gamma}_{n+1}(z)}{{\Gamma}_{n}(z)}, \quad z \in
\dC, \quad n \in \dN, \nonumber \\
{\Gamma}_{1}(z) &=& {{\Gamma}(z)}, \\
{\Gamma}_{n}(1) &=& 1. \nonumber
\end{eqnarray}

\noindent The multiple gamma function, originally introduced over
100 years ago, has significant applications in the connection with
the Riemann Hypothesis. Montgomery \cite{Montgomery} and Sarnak
\cite{Sarnak97} have conjectured that the limiting distribution of
the non-trivial zeros of the Riemann zeta function is the same as
that of the eigenphases of matrices in the CUE (the circular
unitary ensemble). It has been shown in works by Mehta, Sarnak,
Conrey, Keating, and Snaith that a closed representation for
statistical averages over CUE of $N \times N$ unitary matrices,
when $N\rightarrow  \infty$ can be expressed in terms of the
Barnes function $G(z) = 1/{{\Gamma}_2(z)}$, defined by

\begin{eqnarray}
G(z+1) &=& G(z) \, \Gamma(z), \quad z \in
\dC, \\
G(1) &=& 1.\nonumber
\end{eqnarray}

\noindent Keating and Snaith \cite{Keating00,Keating01}
conjectured the following relationship between the moments of $|
\zeta \left(\frac{1}{2}+i \,t\right)|$, averaged over $t$, and
characteristic polynomials, averaged over the CUE:

\begin{equation}
\label{moments}
\begin{array}{c}
\displaystyle \log\bigg(\frac{1}{a(\lambda )}
\lim_{T\rightarrow\infty} \, \frac{1}{{{\log}^{{{\lambda }^2}}}(T)
T}\int _{0}^{T}{{| \zeta
\big(\frac{1}{2}+i t\big) |}^{2\lambda }}\> dt\bigg)=\\
\noalign{\vskip0.5pc} \displaystyle {{\lambda}^2}(\gamma
+1)-2\lambda \sum _{k=2}^{\infty }{{(-\lambda )}^k}
\big({2^k}-1\big) \frac{  \zeta (k)}{ k+1},
\end{array}
\end{equation}
where $a(\lambda )$ is a known function of $\lambda \in \dN$. The
series on the right hand side of (\ref{moments}) is understood in
a sense of analytic continuation (see \cite{Analysis98} for the
method of evaluation of such sums), provided by the Barnes
function

\begin{eqnarray}
\label{Gsum} 2 \sum _{k=2}^{\infty }{{(-z)}^k} \frac{ \zeta
(k)}{k+1}=\frac{2}{z} \log G(z+1) + z(\gamma + 1) -\log(2 \pi ) +
1,\\
\quad |z| <1, \nonumber
\end{eqnarray}

\noindent where $\gamma$ denotes the Euler-Mascheroni constant.
Conversely, we find that the moments of $| \zeta
\left(\frac{1}{2}+i \,t\right)|$ in (\ref{moments}) are just a
ratio of two Barnes functions:

\[
\log\bigg(\frac{1}{a(\lambda )} \lim_{T\rightarrow\infty} \,
\frac{1}{{{\log}^{{{\lambda }^2}}}(T) T}\int _{0}^{T}{{| \zeta
\big(\frac{1}{2}+i \,t\big) |}^{2\lambda }}\> dt\bigg)=
\frac{G(\lambda + 1)^{2}}{G(2 \lambda+1)}
\]
The evidence in support of Keating and Snaith's conjecture is
confirmed by a few particular cases $\lambda =1, 2, 3, 4$, and the
numerical experiment, conducted by Odlyzko \cite{Odlyzko}, for $T$
up to $10^{22}$th zero of the Riemann zeta function.

The identity (\ref{Gsum}) can be further generalized to

\begin{eqnarray}
\label{G3sum} \sum _{k=2}^{\infty }{{(-z)}^k} \frac{ \zeta
(k)}{k+r-1},\quad  r \in \dQ, \quad |z| <1 \nonumber
\end{eqnarray}
and then evaluated in terms of the $r$-tiple gamma functions. For
instance, with $r=3$ we find

\begin{eqnarray}
\label{Gamma3sum} \sum _{k=2}^{\infty }{{(-z)}^k} \frac{ \zeta
(k)}{k+2} &=& \frac{2}{z^2} \log \Gamma_{3}(z+1) + \frac{1}{z^2}
\log G(z+1) + \\
\frac{6 z^{2} + 3 z -1}{12 z} &+& \frac{\gamma z}{3} - \frac{\log
(2 \pi)}{2} - \frac{2 \,\zeta^{\prime}(-1)}{z}, \quad |z| <1.
\nonumber
\end{eqnarray}

\medskip
    In this paper we aim at developing a mathematical foundation for
symbolic computation of special classes of infinite series and
products. Generally speaking, all series (subject to convergence)
of the form:

\[\sum _{j=1}^{\infty }R(j)\, {\log ^m}P(j), \quad m \in \dN, \]
where $R(z)$ and $P(z)$ are polynomials, can be expressed in a
closed form by means of the multiple gamma function, which may
further simplify to elementary functions. This algorithm also
complements the work previously started by Adamchik and Srivastava
\cite{Analysis98} for symbolic summation of series involving the
Riemann zeta function. The algorithm will expand the capabilities
of existing software packages for non-hypergeometric summation.

\section{Asymptotics of ${\zeta}^{\prime}(-\lambda ,z)$}

In this section, based on the integral (\ref{Hermite}), we derive
the asymptotic expansion of ${\zeta}^{\prime}(t,z) = \frac{d}{d
t}\,\zeta (t,z)$, when $t=-\lambda, \> \lambda \in \dN_{0}$ and
$z\rightarrow \infty$. Differentiating both sides of
(\ref{Hermite}) with respect to $s$, we obtain

\begin{eqnarray}
\label{int0} \zeta^{\prime}(-\lambda ,z) &=& \frac{{z^{\lambda
+1}}}{\lambda +1}\log  z-\frac{{z^{\lambda }}}{2} \log
z-\frac{{z^{\lambda +1}}}{{{(\lambda +1)}^2}}\nonumber\\
 &+& 2
\int _{0}^{\infty }\frac{{{\tan}^{-1}}\big(\frac{x}{z}\big)
\cos\big(\lambda \,
{{\tan}^{-1}}\big(\frac{x}{z}\big)\big)}{({e^{2 \pi
x}}-1){{({x^2}+{z^2})}^{-\lambda /2}} }\> dx \\ \nonumber
 &+& \int
_{0}^{\infty }\frac{ \log({x^2}+{z^2}) \sin\big(\lambda \,
{{\tan}^{-1}}\big(\frac{x}{z}\big)\big)}{({e^{2 \pi
x}}-1){{({x^2}+{z^2})}^{-\lambda /2}} }\> dx
\end{eqnarray}
Next, we expand the integrands into the Taylor series with respect
to $z$. Taking into account that

\[
{{\big({x^2}+{z^2}\big)}^{\lambda /2}} \cos\big(\lambda
\,{{\tan}^{-1}}\big(\frac{x}{z}\big)\big)=\sum _{k=0}^{\lambda /2}
{{(-1)}^k} {\lambda \choose 2 k}{z^{\lambda -2 k}} {x^{2 k}}
\]
and

\[
{{\big({x^2}+{z^2}\big)}^{\lambda /2}} \sin\big(\lambda\,
{{\tan}^{-1}}\big(\frac{x}{z}\big)\big)=\sum _{k=0}^{\lambda /2}
{{(-1)}^k}{\lambda \choose 2 k+1}{z^{\lambda-2 k -1}} {x^{2 k+1}}
\]

\noindent we compute

\begin{equation}
\label{int1}
\begin{array}{c}
\displaystyle \int _{0}^{\infty
}\frac{{{\tan}^{-1}}\big(\frac{x}{z}\big)\, \cos\big(\lambda\,
{{\tan}^{-1}}\big(\frac{x}{z}\big)\big)}{({e^{2
\pi  x}}-1){{({x^2}+{z^2})}^{-\lambda /2}} }\> dx=   \\
\noalign{\vskip 1.0pc} \displaystyle
      \sum_{k=0}^{\lambda /2}{z^{\lambda -2 k}} {{(-1)}^k}
{\lambda \choose 2 k} \int _{0}^{\infty }\frac{{x^{2 k}}\,
{{\tan}^{-1}}\big(\frac{x}{z}\big)}{{e^{2 \pi  x}}-1}\> dx
\end{array}
\end{equation}
and

\begin{equation}
\label{int2}
\begin{array}{c}
\displaystyle \int _{0}^{\infty }\frac{ \log({x^2}+{z^2})
\sin\big(\lambda
\,{{\tan}^{-1}}\big(\frac{x}{z}\big)\big)}{({e^{2
\pi  x}}-1){{({x^2}+{z^2})}^{-\lambda /2}} }\> dx=    \\

\noalign{\vskip 1.0pc} \displaystyle

-\sum _{k=0}^{\lambda /2}{z^{-2 k+\lambda -1}} {{(-1)}^k}{\lambda
\choose 2k+1} \int _{0}^{\infty }\frac{{x^{2 k+1}}
\log({x^2}+{z^2})}{{e^{2 \pi  x}}-1}\> dx.
\end{array}
\end{equation}

\noindent In the next step, we find asymptotic expansions of
integrals on the right hand side of (\ref{int1}) and (\ref{int2})
when $z\rightarrow \infty$. Expanding
${\tan}^{-1}\big(\frac{x}{z}\big)$ and $\log\big({x^2}+{z^2}\big)$
into the Taylor series with respect to $x$ and performing term by
term integration, we arrive at

\[
\int _{0}^{\infty }\frac{{x^{2 k}}\,
{{\tan}^{-1}}\big(\frac{x}{z}\big)}{{e^{2 \pi x}}-1}\>d x=\sum
_{j=0}^{r}\frac{{{(-1)}^k} {z^{-2 j-1}} {B_{2 (j+k+1)}}}{4 (2 j+1)
(j+k+1)}+ O\big(\frac{1}{{z^{ 2r+3}}}\big),
\]

\[
\int _{0}^{\infty }\frac{{x^{2 k+1}} \log(1+{x^2}/{z^2})}{{e^{2
\pi  x}}-1}\>d x=-\sum _{j=1}^{r}\frac{{{(-1)}^k} {z^{-2 j}} {B_{2
(j+k+1)}}}{4 j (j+k+1)}+ O\big(\frac{1}{{z^{ 2r+2}}}\big),
\]
where $B_{k}$ are the Bernoulli numbers. Substituting these into
(\ref{int1}) and (\ref{int2}) and the latters into (\ref{int0}),
after some tedious algebra, we obtain

\begin{prop}
The derivatives of the Hurwitz zeta function have the following
asymptotic expansion when $z\rightarrow \infty$:

\begin{eqnarray}
\label{prop1}
\zeta ^{\prime}(-\lambda ,z) &=& \frac{{z^{\lambda +1}}}{\lambda
+1}\log  z-\frac{{z^{\lambda }}}{2} \log
z- \frac{{z^{\lambda+1}}}{{{(\lambda +1)}^2}}   \nonumber \\
&+&\sum _{j=1}^{r}\frac{{z^{-2
j+\lambda +1}}}{(2 j)!}{B_{2 j}}\sum _{k=0}^{2 j-1}(k+1)
{\left[{{2j-1} \atop {k+1}}\right]}{{(-\lambda )}^k}   \\\nonumber
\displaystyle
&-&\log  z \sum _{j=1}^{r}\frac{{z^{-2 j+\lambda +1}} {B_{2 j}}
{{(-\lambda )}_{2 j-1}}}{(2 j)!} + O\big(\frac{1}{{z^{ 2\,r-\lambda
+1}}}\big)
\end{eqnarray}

\noindent where ${\zeta}^{\prime}(t,z) = \frac{d}{d t}\,\zeta
(t,z)$, $(-\lambda )_{k} = (-\lambda )(-\lambda +1)\ldots(-\lambda
+k-1)$ is the Pochhammer symbol, and $\left[{n \atop k}\right]$
are the Stirling cycle numbers, defined recursively \cite{Rosen} by

\begin{equation}
\label{cycle} {\left[{n \atop k}\right]} = (n-1) \, {\left[{n-1
\atop k}\right]} + {\left[{n-1 \atop k-1}\right]} , \quad
{\left[{n \atop 0}\right]} = {\Bigg\{\begin{array}[c]{c}
    1, \quad n = 0, \\ \noalign{\vskip0.5pc}
    0, \quad n \neq 0
  \end{array}}
\end{equation}

\end{prop}

\noindent
For $\lambda = 0$ and $\lambda = 1$, formula (\ref{prop1}) yields the well-know
asymptotics \cite{Magnus}:

\begin{equation}
\label{asymp}
\begin{array}{l}
\displaystyle
\zeta^{ \prime}(0,N) = N \log N - \frac{\log N}{2} -N+
O\big(\frac{1}{{N}}\big), \\
\noalign{\vskip1.5pc} \displaystyle
\zeta^{ \prime}(-1,N) = \left(\frac{{N^2}}{2}-\frac{N}{2}+\frac{1}{12} \right)\log N
-\frac{{N^2}}{4}+\frac{1}{12}+O\big(\frac{1}{{N^2}}\big)
\end{array}
\end{equation}

\noindent
In the similar way, we can derive the asymptotic expansions for higher order derivatives
$\zeta^{(m)}(-\lambda, N)$, where $N \rightarrow \infty$ and $\lambda \in \dN^{+}$.

\section{Multiple gamma, zeta and the Hurwitz functions}
In [20], Vardi expressed the ${\Gamma}_n$ function in terms of the
multiple zeta function ${\zeta}_n(s,z)$, as

\begin{equation}
\label{logGn}
(-1)^{n} \log {\Gamma_n}(z)=-{\lim_{s\rightarrow 0}}
\left(\frac{\partial \,{{\zeta }_n}(s,z)}{\partial s}\right)-\sum
_{k=1}^{n}{{(-1)}^k} {z \choose {k-1}} {R_{n+1-k}}
\end{equation}

where

\begin{equation}
\label{Rn} {R_n}=\sum  _{k=1 }^{n }\lim_{s\rightarrow 0}
\left(\frac{\partial \,{{\zeta }_k}(s,1)}{\partial s}\right)
\end{equation}

and

\begin{eqnarray}
\label{multipleZ}
{{\zeta }_n}(s,z) &=& \sum _{{k_1}=0}^{\infty }\ldots
\sum _{{k_n}=0}^{\infty
}\frac{1}{{{({k_1}+{k_2}+...+{k_n}+z)}^s}} \nonumber\\
 &=&
\sum _{k=0}^{\infty
}\frac{1}{{{(k+z)}^s}}\,{k+n-1 \choose n-1}
\end{eqnarray}

The aim of this section is to find a closed form representation
for $\log {\Gamma}_n(z)$ in finite terms of the Hurwitz zeta
function, and vice versa.

\begin{prop}
The multiple gamma function ${\Gamma}_{n}(z)$ may be expressed by
means of the derivatives of  the Hurwitz zeta function:

\begin{equation}
\label{prop2}  \log {{\Gamma}_n}(z)=\frac{(-1)^{n}}{(n-1)!} \sum
_{k=0}^{n-1}{P_{k,n}}(z)
\bigg(\zeta^{\prime}(-k)-\zeta^{\prime}(-k, z) \bigg), \quad
\Re(z)
> 0,
\end{equation}
where the polynomials $P_{k, n}(z)$  are defined by

\begin{equation}
\label{P} {P_{k, n}}(z)=\sum _{j=k+1}^{n}{{(-z)}^{j-k-1}} {j-1
\choose k} {\left[ n \atop j\right]}
\end{equation}
where ${\left[ n \atop j\right]}$ are the Stirling cycle numbers.

\end{prop}

The polynomials ${P_{k,n}}(z)$ can be envisaged as the generalized
Stirling polynomials of the first kind, generated by

\[\prod _{k=1}^{n-1}(k+x-z)=\sum _{k=0}^{n-1}{P_{k,
n}}(z){x^k}\]
with the following two alternative forms of representation

\begin{equation}
\label{P1} {P_{k, n}}(z)=\frac{{{(-1)}^k}}{k!} \frac{{{\partial
}^{n-1}}}{\partial {y^{n-1}}}\frac{{{\log}^k}(1-y)}{
{{(1-y)}^{1-z}}}\bigg|_{y\rightarrow 0}
\end{equation}
and

\begin{equation}
\label{P2} {P_{k, n}}(z)=\sum_{i=k+1}^{n}{z \choose n-i}\frac{
(n-1)! }{(i-1)!} {\left[ {i \atop k+1} \right]}
\end{equation}
For $z=1$ the polynomials ${P_{k,n}}(z)$ are simplified to the
Stirling numbers:

\begin{equation}
\label{P_at1}
{P_{k, n}}(1)= {\left[ {n-1 \atop k} \right]}
\end{equation}

\noindent
Interestingly, the polynomials ${P_{k,n}}(z)$ were first considered by
Mitrinovic \cite{Mitrinovic} 40 years ago (with no relation to the multiple
gamma function) as a possible generalization of the Stirling
numbers.

\medskip
    The proof of Proposition 2 is based on the following two lemmas.

\begin{lem}
The multiple zeta function ${{\zeta }_n}(s,z)$ defined by
(\ref{multipleZ}) may be expressed by means of the Hurwitz function

\begin{eqnarray}
\label{lem1} {{\zeta }_n}(s,z)=\frac{1}{(n-1)!}\sum
_{j=0}^{n-1}{P_{j, n}}(z)\, \zeta (s-j,z)
\end{eqnarray}
\end{lem}

\noindent
\begin{pf*}{Proof}
Recall the definition of the Stirling cycle numbers

\[
{k+n-1 \choose n-1} = \frac{1}{(n-1)!}\sum _{i=0}^{n}{k^{i-1}}
{\left[ n \atop i \right]}
\]
\noindent
Expanding ${k^{i-1}}={{((k+z)-z)}^{i-1}}$ by the binomial theorem,
implies that

\[
{k+n-1 \choose n-1} = \frac{1}{(n-1)!} \sum_{i=0}^{n}\sum
_{j=0}^{i-1}{i-1 \choose j} (k+z)^{j} {{(-z)}^{i-j-1}} {\left[ n
\atop i \right]}
\]
\noindent
Interchanging the order of summation and making use of (\ref{P}),
leads us to

\[
{k+n-1 \choose n-1} = \frac{1}{(n-1)!}\sum _{j=0}^{n} {P_{j,
n}}(z) {{(k+z)}^j}
\]
\noindent
Finally, we substitute this binomial coefficients representation
into (\ref{multipleZ}):

\begin{eqnarray*}
{{\zeta }_n}(s,z) &=& \sum _{k=0}^{\infty }\frac{1}{{{(k+z)}^s}}
{k+n-1 \choose n-1} \\
&=& \frac{1}{(n-1)!}\sum _{j=0}^{n}{P_{j,
n}}(z)\sum _{k=0}^{\infty }\frac{1}{{{(k+z)}^{s-j}}}
\end{eqnarray*}
\noindent
We complete the proof by evaluating the inner sum to the Hurwitz
zeta function.\qed
\end{pf*}

\begin{lem}
The function $R_n$ defined by (\ref{Rn})  may be expressed by
means of the derivatives of the Riemann zeta function

\begin{equation}
\label{lem2} {R_n}=\frac{1}{(n-1)!}\sum _{k=0}^{n-1}
\zeta^{\prime}(-k) {\left[ {n \atop k+1} \right]}
\end{equation}
\end{lem}

\noindent
\begin{pf*}{Proof}
\noindent Using Lemma 1 with $z=1$ and formula (\ref{P_at1}), we
get

\[
{{\zeta }_k}(s,1)= \frac{1}{(k-1)!}\sum _{j=0}^{k} {\left[ {k-1 \atop j}
\right]} \zeta (s-j)
\]
which, upon differentiation, yields

\begin{eqnarray}
\label{lim_lem2}
\lim_{s \rightarrow 0}
 \frac{\partial{{\zeta }_k}(s,1)}{\partial
s} = \zeta^{\prime}(0)\, \delta_{0,k-1} +
\frac{1}{(k-1)!}\sum_{j=1}^{k}{\left[ {k-1 \atop j} \right]}
\zeta^{\prime}(-j)
\end{eqnarray}
\noindent
Here ${\delta }_{k,n}$ is the Kronecker delta function. Next, we
sum both parts of identity (\ref{lim_lem2}) with respect to $k$:

\[
{R_n}=\sum _{k=1}^{n}\lim_{s \rightarrow 0} \frac{\partial {{\zeta
}_k}(s,1)}{\partial s} = \zeta^{\prime}(0) + \sum _{j=1 }^{n
}\zeta^{\prime}(-j)\sum _{k=j}^{n}\frac{1}{(k-1)!} {\left[ {k-1
\atop j} \right]}
\]
\noindent
Observing, that the inner sum can be expressed via the Stirling
numbers:

\[
\sum _{k=j}^{n}\frac{1}{(k-1)!}
{\left[ {k-1 \atop j} \right]} =\frac{1}{(n-1)!} {\left[ {n \atop
j+1} \right]}
\]
we complete the proof.\qed
\end{pf*}

\begin{pf*}{Proof of Proposition 2}
The derivation of (\ref{prop2}) immediately follows from formula
(\ref{logGn}), by applying Lemma 1 and Lemma 2:

\[
\begin{array}[l]{c}
\displaystyle

\log {\Gamma_n}(z)=-\frac{(-1)^{n} }{(n-1)!}\sum _{j=0}^{n-1}{P_{j, n}}(z)\,
\zeta^{\prime}(-j,z)\,+  \\
\noalign{\vskip0.5pc} \displaystyle

\sum _{j=0}^{n-1}{{(-1)}^j} {z \choose j}
 \frac{1}{(n-j-1)!}\sum _{k=0}^{n-j-1} \zeta^{\prime}(-k)
 {\left[ {n-j \atop k+1} \right]}
\end{array}
\]

\noindent
Interchanging the order of summation, we obtain

\[
\begin{array}[l]{c}
\displaystyle

\log {\Gamma}_{n}(z)=\frac{(-1)^{n}}{(n-1)!}\sum
_{j=0}^{n-1} \Bigg(-P_{j,n}(z) \zeta^{\prime}(-j,z)\,+  \\
\noalign{\vskip0.5pc} \displaystyle

\zeta^{\prime}(-k)\sum _{j=0}^{n-k} {{(-1)}^j} {z \choose j}
 \frac{(n-1)!}{(n-j-1)!}
 {\left[ {n-j \atop k+1} \right]} \Bigg)
\end{array}
\]
\noindent
Finally, we note that the above inner sum is equal to
representation (\ref{P2}) with a reversed order of summation. This
completes the proof.\qed
\end{pf*}

    In the second part of this section we invert Proposition 2 and
express the derivatives of the Hurwitz function in finite terms of
the multiple gamma function. The problem of representing
derivatives in terms of other special functions was originally
addressed by Ramanujan who studied functional and asymptotic
properties of this function  \cite{Berndt1}

\[
\phi_{k}(z) = {{\zeta }^{\prime }}(-k,z+1)-{{\zeta }^{\prime
}}(-k).
\]
\noindent
It is unknown if Ramanujan was aware of Barnes' results regarding
the multiple gamma function.

\begin{prop}
The derivatives of the Hurwitz zeta function may be expressed by
means of the multiple gamma function ${\Gamma}_{n}(z)$

\begin{equation}
\label{prop3} \zeta^{\prime}(-n,z) - \zeta^{\prime}(-n) = (-1)^{n}
\sum _{k=0}^{n}k!\>{Q_{ k, n}}(z) \log {{\Gamma}_{k+1}}(z), \quad
\Re(z) > 0
\end{equation}
\noindent
where the polynomials ${Q_{ k, n}}(z)$ are defined by

\begin{equation}
\label{Q} Q_{ k, n}(z)=\sum _{j=k}^{n}(1-z)^{n-j} {\left(n \atop j
\right)} {\left\{ {j \atop k} \right \}}
\end{equation}
\noindent
and ${\left\{ {j \atop k} \right \}}$  are the Stirling subset
numbers, defined by \cite{Rosen}

\begin{equation}
\label{subset} {\left\{{n \atop k}\right\}} = k \, {\left\{{n-1
\atop k}\right\}} + {\left\{{n-1 \atop k-1}\right\}} , \quad
{\left\{{n \atop 0}\right\}} = {\Bigg\{\begin{array}[c]{c}
    1, \quad n = 0, \\\noalign{\vskip0.5pc}
    0, \quad n \neq 0
  \end{array}}
\end{equation}

\end{prop}

The immediate particular case of (\ref{prop3}) is the formula for
the triple gamma function:

\begin{eqnarray}
\zeta^{\prime}(-2,z) - \zeta^{\prime}(-2) &=& 2 \log \Gamma_{3}(z)
+ (3-2z) \log G(z) \\
\noalign{\vskip0.5pc} &+& (1-z)^{2} \log \Gamma(z), \quad \Re(z)
> 0 \nonumber
\end{eqnarray}

For $z=0$ the polynomials ${Q_{k,n}}(z)$ simplify to the Stirling
subset numbers:

\[
{Q_{ k, n}}(0)  = \sum_{j=k}^{n}{n \choose j} {\left\{ {j \atop k}
\right\}} = {\left\{ {n+1 \atop k+1} \right\}}
\]

\begin{lem}
Polynomials  ${P_{j,n}}(z)$ and ${Q_{k,j}}(z)$ satisfy the
following discrete orthogonality relation:

\begin{equation}
\label{lem3} \sum _{j=k}^{n-1}{{(-1)}^{j-k}}{Q_{ k, j}}(z)
 {P_{j,n}}(z) = {{\delta }_{k,n-1}}
 \end{equation}

\noindent where ${{\delta }_{k,n-1}}$ is the Kronecker delta.
\end{lem}

\medskip
\begin{pf*}{Proof} The proof is primarily based on the discrete orthogonality
relation for the Stirling numbers:

\[
\sum _{j=0}^{n} (-1)^{m+j} {\left  [ {j \atop m} \right ]} {\left
\{ {n \atop j} \right \}} = \delta_{m, n}
\]
\qed
\end{pf*}

\begin{pf*}{Proof of Proposition 3}
The proof readily follows from (\ref{prop3}) by replacing $\log
{\Gamma}_{k+1}(z)$ by (\ref{prop2}) and then using lemma 3.\qed
\end{pf*}

Polynomials ${Q_{k,n}}(z)$ generate a great deal of new integer
sequences as well as closed form representations for the existing
ones.  Here are two examples from the online encyclopedia of
integer sequences \cite{Sloan}.

\bigskip

{\bf Example 1}

The sequence A021424: $1, 16, 170, 1520, 12411, 96096$ made of
coefficients by $x$ in the expansion of $1/((1-x) (1-3 x) (1-5
x) (1-7 x))$. This sequence has a closed form representation

\[
{2^{n-3}} {Q_{3, n}}\big(\frac{1}{2}\big)=\frac{1}{48}
\big({7^n}-3 \, {5^n}+{3^{n+1}}-1\big)
\]

\medskip

{\bf Example 2}

Riordan's sequence A000554 (labeled trees of diameter 3 with
$n$ nodes)$: 12, 60, 210, 630, 1736, 4536, 11430$. This sequence
is generated by the coefficients of $z^2$ in the polynomial
$2*{Q_{2, n}}(1-z)$, which is also can be written as

\[
\frac{{{\partial }^2}{Q_{2, n}}(1-z)}{\partial {z^2}}\bigg|_{z
\rightarrow 0} =2 {n \choose n-2}{\left\{ {n-2 \atop 2}\right \}}
\]

\section{Summation}

Consider the class of infinite sums

\begin{equation}
\label{genser} \sum _{k=1}^{\infty }R(k)\log  P(k)
\end{equation}
\noindent
where $R(z)$ and $P(z)$ are polynomials. We will show that all
such sums can be expressed in finite terms of the multiple gamma
function. Clearly, (\ref{genser}) is a linear combination (subject
to the branch cut of log) of

\begin{equation}
\label{genser1} {\Phi_N}(z)=\lim_{N \rightarrow \infty}
\sum_{k=1}^{N}{k^p} \log (k+z)
\end{equation}
\noindent
Thus, it is sufficient to find a closed form for (\ref{genser1}).
Expanding

\[{k^p}={{((k+x)-x)}^p}\]
by virtue of the binomial theorem, we get

\[{{\Phi }_N}(z)=\lim_{N \rightarrow \infty}
\sum_{j=0}^{p} {p \choose j} (-1)^j z^{p-j} {F_N}(z) \]
\noindent
where

\begin{equation}
\label{fn} {F_N}(z)=\sum _{k=0}^{N}{{(k+z)}^j}\log
(k+z)\end{equation}

Assuming the analytic property of the Hurwitz zeta function,
(\ref{fn}) can be written in the form

\[{F_N}(z)=-\lim_{s\rightarrow \infty}
\frac{d}{ds}\sum _{k=0}^{N}\frac{1}{(k+z)^s} = \zeta^{\prime}(
-j,N+z+1) - \zeta^{\prime}(-j,z) \]
\noindent
and therefore

\begin{equation}
\label{regular}
\begin{array}{l}
\displaystyle
 {{\Phi}_N}(z)=\lim_{N \rightarrow \infty}
\sum_{j=0}^{p} {p \choose j} (-1)^{j+1} z^{p-j}  \\
\noalign{\vskip0.5pc}{\hskip3.5pc}  \displaystyle

\Bigg(\zeta^{\prime}(-j,N+z+1)- \zeta^{\prime}(-j,z)\Bigg)
\end{array}
\end{equation}
\noindent
Further, we replace the first cumulant $\zeta^{\prime}(-j,z)$ by a
linear combination of the multiple gamma function (see Proposition
3), and the second cumulant $\zeta^{\prime}(-j,N+z+1)$ by
asymptotic (\ref{prop1}). When the original series (\ref{genser})
converges, the $N$-dependent term will vanish as $N\rightarrow
\infty$. However, if series (\ref{genser}) diverges, formula
(\ref{regular}) combined with (\ref{prop1}) will provide a zeta
regularization of (\ref{genser}) in the Hadamard sense.

\medskip

In order to make this idea clear, we will consider a few examples.

\subsection{Dilcher's sum}

In \cite{Dilcher94}, Dilcher introduced a particular generalization of the
Euler gamma function $\Gamma(x)$ which is related to the Stieltjes
constants $\gamma_{k}$, as $\Gamma(x)$ is related to the Euler
$\gamma$ constant. The Stieltjes constants $\gamma_{k}$ are
defined as coefficients in the Laurent expansion of the Riemann
zeta function $\zeta (s)$ at the simple pole $s=1$:

\[
\zeta (s)=\frac{1}{s-1}+\sum
_{k=0}^{\infty}\frac{{{(-1)}^k}}{k!}{{\gamma }_k} {{(s-1)}^k} \]
\noindent
The Euler $\gamma$ constant is a particular case of the Stieltjes
constants $\gamma_{0} = \gamma$. Dilcher derived a few basic
properties of the generalized gamma function as well as its
asymptotic expansion when $z\rightarrow \infty$. The following
infinite series appears as a constant term in the asymptotic
expansion:

\begin{equation}
\label{Dk}
{D_k}=\sum _{j=1}^{\infty
}{{\log}^k}\big(j+\frac{1}{2}\big)-2 {\log
^k}j+{{\log}^k}\big(j-\frac{1}{2}\big)
\end{equation}

\noindent
For  $k=1$ the sum (\ref{Dk}) simplifies to ${D_1}=\log\left(2/\pi\right)$. For other values of $k>1$, the closed form was unknown. In this
section, we apply the technique of the multiple gamma function and
evaluate the logarithmic series (\ref{Dk}) in closed form.

\medskip
We begin with a finite sum

\begin{equation}
\label{limDk} {D_k}=\lim_{N\rightarrow \infty}\sum
_{j=1}^{N}{{\log}^k}\big(j+\frac{1}{2}\big)-2 \,{\log
^k}j+{{\log}^k}\big(j-\frac{1}{2}\big)
\end{equation}
\noindent
and then rewrite it in terms of the Hurwitz function. The
following chain of operations is valid (assuming the analytic
property of the Hurwitz function)

\begin{equation}
\label{tozeta}
\begin{array}[l]{c}
\displaystyle
(-1)^k \sum_{j=1}^{N} \log^k (j+c) =
\lim_{s\rightarrow 0} \frac{{{\partial }^k}}{\partial {s^k}}

\sum_{j=1}^{N} \frac{1}{{{(j+c)}^s}} = \\
\noalign{\vskip0.5pc} \displaystyle

\lim_{s\rightarrow 0} \frac{{{\partial }^k}}{\partial {s^k}}
\Bigg(\sum_{j=1}^{\infty }\frac{1}{{{(j+c)}^s}}-\sum
_{j=N+1}^{\infty}\frac{1}{{{(j+c)}^s}}\Bigg)  = \\
\noalign{\vskip0.5pc} \displaystyle

\zeta^{k}(0,c+1)-\zeta^{k}(0,c+N+1)
\end{array}
\end{equation}

\noindent
Applying this to (\ref{limDk}), yields

\begin{equation}
\label{limDk1}
\begin{array}{l}
\displaystyle

{{(-1)}^k}{D_k}=-2 \zeta^{k}(0,1) + \zeta^{k}(0,\frac{1}{2}) +
\zeta^{k}(0,\frac{3}{2}) \,+\\
\noalign{\vskip0.5pc} \displaystyle \hspace{2.5pc}

 \lim_{N\rightarrow
\infty} \left(2\,\zeta^{k}(0,N+1) - \zeta^{k}(0,N+\frac{1}{2}) -
\zeta^{k}(0,N+\frac{3}{2}) \right)
\end{array}
\end{equation}
\noindent
Now we use the asymptotic expansion of $\zeta^{k}(0,N)$, $N \rightarrow \infty$. From the
Hermite integral (\ref{Hermite}) it follows that the dominant
asymptotic term of $\zeta^{k}(0,N)$ comes from the first two terms

\[
\zeta(s, z) = \frac{z^{-s}} {2} + \frac{z^{1-s}} {s-1} + \ldots
\]
\noindent
This leads to the following asymptotic expansion $N \rightarrow
\infty$

\begin{equation}
\begin{array}{l}
\displaystyle
\zeta^{k}(0,N) =N\, \sum _{j=0}^{k}{{(-1)}^j}
{{(k-j+1)}_j} \, {\log ^{k-j+1}} N +  \\
\noalign{\vskip0.5pc} \displaystyle \hspace{4.5pc}

\frac{{{(-1)}^k}}{2} \,{\log ^k} N +O\big(\frac{1}{N}\big),
\end{array}
\end{equation}
\noindent
where $(k-j+1)_j$ is the Pochhammer symbol. Thus, to prove that
the limiting part in (\ref{limDk1}) vanishes to zero, it is
sufficient to show that

\begin{equation}
\label{lim_0}
\begin{array}{c}
\displaystyle

\lim_{N \rightarrow \infty} \Bigg[2(N+1){\log ^{j}} (N+1) \,- \hspace{4.5pc}\\
\noalign{\vskip0.5pc} \displaystyle

\left(N+\frac{1}{2}\right) {\log ^{j}}
\left(N+\frac{1}{2}\right) \,-
\left(N+\frac{3}{2}\right){\log^{j}}
\left(N+\frac{3}{2}\right)\Bigg]=0
\end{array}
\end{equation}

\noindent
for any $j = 1, 2,  \ldots, k+1$. This follows immediately from the
fact that

\[
\log(N+c)=\log N +\frac{c}{N} + O\left(\frac{1}{{N^2}}\right),
\quad N\rightarrow \infty
\]
\noindent
and the sum of coefficients by $\log^{k-j+1}(N)$ in (\ref{lim_0})
is zero.

Therefore, we have

\begin{equation}
\label{Dkfinal}
D_{k} =(-1)^{k} \left(-\log^{k} 2 - 2 \,\zeta^{k}(0)
+ 2 \,\zeta^{k}(0,\frac{1}{2}) \right)
\end{equation}
\noindent
Note, $\zeta^{k}(0,\frac{1}{2})$ in the right hand-side of
(\ref{Dkfinal}) is a linear combination of $\zeta^{k}(0)$:

\[
\zeta^{k}(0,\frac{1}{2}) = -\frac{1}{2}{\log ^k}2 +
\sum_{j=1}^{k-1} {k \choose j} \zeta^{k-j}(0) \log^{j} 2
\]
\noindent
Combining this with (\ref{Dkfinal}), yields

\begin{prop}
The infinite series (\ref{Dk}) can be expressed in finite terms of
the derivatives of the Riemann zeta function

\begin{equation}
\label{prop4}
\frac{(-1)^{k}}{2} D_{k} = -\log^{k} 2 - \zeta^{k}(0) +
\sum_{j=1}^{k-1} {k \choose j} \zeta^{k-j}(0) \log^{j} 2
\end{equation}
\end{prop}

\noindent
In particular, for $k=2$ and $k=3$ we have

\begin{eqnarray}
\label{D2}
{D_2}&=&\sum _{j=1}^{\infty }{{\log}^2}\left(j+\frac{1}{2}\right)-2\,
{\log ^2}j+{{\log}^2}\left(j-\frac{1}{2}\right) \nonumber\\
&=&\frac{{{\pi }^2}}{12}+{\log ^2}\pi -3 \,{\log ^2}2-{{\gamma }^2}-2
{{\gamma }_1}
\end{eqnarray}

\begin{eqnarray}
\label{D3}
{D_3}= 3 D_{2} \log 2 + 2 \, \zeta^{\prime\prime\prime}(0) + 9 \log(2 \pi) \log^{2}2 + \log^{3} 4
\end{eqnarray}
respectively.

\medskip
Note that by the functional equation
\[
\zeta (s)={2^s} {{\pi }^{s-1}} \Gamma (1-s) \sin\big(\frac{\pi
s}{2}\big) \, \zeta (1-s), \quad \Re(s)\neq1
\]
the Riemann zeta function admits a meromorphic continuation to the entire complex plane. Differentiating this with respect to $s$ and setting $s=0$, allows us to express derivatives $\zeta^{k}(0)$ in finite terms of the Stieltjes constants (see \cite{Apostol85} for a closed form evaluation). Therefore, formulas like (\ref{D2}) and (\ref{D3}) can be used for numeric computation of $\zeta^{k}(0)$ as well as the Stieljes constants $\gamma_{k}$.

\subsection{Melzak's Product}
Let us consider the infinite product

\begin{equation}
\label{Melzak} \lim_{N \rightarrow \infty} \prod_{k=1}^{2 N}
{{\left(1+\frac{2\,x}{k}\right)}^{-k (-1)^k}}
\end{equation}

\noindent
and evaluate it in closed form. Originally, this product was computed by Melzak \cite{Melzak} in a particular case when $x = 1$

\begin{equation}
\label{Melzak1} \lim_{N \rightarrow \infty} \prod_{k=1}^{2 N}
{{\left(1+\frac{2}{k}\right)}^{-k (-1)^k}} = \frac{\pi}{2\,e}
\end{equation}
Later on, P. Borwein and W. Dykshoorn \cite{Dykshoorn} generalized
(\ref{Melzak1}) and computed it in terms of Bendersky's function,
defined by (see [26])

\[
1^{1} 2^{2} \ldots n^{n} = \Gamma_{1}(n+1),
\]
which is closely related to the Barnes function $G$

\[
\Gamma_{1}(n+1) = \frac{n!^{n}}{G(n+1)}
\]
Note that (\ref{Melzak}) can be written in an alternative form
by changing the product limit from $2N$ to $2N+1$:

\[
\lim_{N \rightarrow \infty} \prod_{k=1}^{2 N+1}
{{\left(1+\frac{2\,x}{k}\right)}^{-k (-1)^k}} = e^{2 \,x} \lim_{N \rightarrow \infty} \prod_{k=1}^{2 N}
{{\left(1+\frac{2\,x}{k}\right)}^{-k (-1)^k}}
\]
\noindent
In this section, using the multiple gamma function technique, we
reevaluate product (\ref{Melzak}), derive its few particular cases in terms of known constants and then generalize it to

\[
\lim_{N \rightarrow \infty} \prod_{k=1}^{2 N}
{{\left(1-\frac{4 x^2}{k^2}\right)}^{-k^2 (-1)^k}}
\]

\begin{prop}
The following identity holds for $\Re(x)>-\frac{1}{2}$

\begin{eqnarray}
\label{prop5} \lim_{N \rightarrow \infty} \prod_{k=1}^{2 N}
{{\left(1+\frac{2x}{k}\right)}^{-k (-1)^k}} = \frac{{e^{-x}}
\Gamma (x+\frac{1}{2}) }{ \Gamma (\frac{1}{2})}
{{\left(\frac{G (x+\frac{1}{2})}{G(x+1) \, G(\frac{1}{2})}\right)}^2}
\end{eqnarray}
\end{prop}

\medskip
\begin{pf*}{Proof} We first convert this product into a finite sum
by applying the logarithm to it:

\begin{equation}
\label{prodtosum}
\begin{array}{l}
\displaystyle
\lim_{N \rightarrow \infty} \prod_{k=1}^{2 N}
{{\left(1+\frac{2x}{k}\right)}^{-k (-1)^k}} = \\
\noalign{\vskip0.5pc} \displaystyle {\hskip6.0pc}
\exp \bigg(\lim_{N \rightarrow \infty}
\sum_{k=1}^{2 N}{(-1)^k} k \,\Big(\log k - \log(k+2x)\Big)\bigg)
\end{array}
\end{equation}
In the next step we split the finite sum in the right hand-side of (\ref{prodtosum}) into three sums:

\[
\begin{array}{r}
\displaystyle
\sum _{k=1}^{2 N}{{(-1)}^k} k \Big(\log  k- \log(k+2x)\Big) = 2 \,x\sum_{k=1}^{2 N}{{(-1)}^k}  \log(k+2 x) \,+\\
\noalign{\vskip0.5pc} \displaystyle

\sum_{k=1}^{2 N}{(-1)^k} k \log k -
\sum_{k=1}^{2 N}{(-1)^k} (k+2 x)\log(k+2x)
\end{array}
\]

\noindent
and evaluate each of them in terms of the Hurwitz function. For
arbitrary $c$ such that $\Re(c)>0$, we find

\[
\begin{array}{r}
\displaystyle
\sum_{k=1}^{2N}{{(-1)}^k} (k+c) \log(k+c)=
-\lim_{s \rightarrow -1}
\frac{\partial }{\partial s}\sum _{k=1}^{2N}
\frac{{{(-1)}^k}}{{{(k+c)}^s}} = \\
\noalign{\vskip0.5pc} \displaystyle

N \log  2 + 2 \,\zeta^{\prime}\left(-1, \frac{c+1}{2} \right) -
2 \,\zeta^{\prime}\left(-1, \frac{c+2}{2} \right) + \\
\noalign{\vskip0.5pc} \displaystyle
2\, \zeta^{\prime}\left(-1, \frac{c+2}{2} +N \right) -
2 \,\zeta^{\prime}\left(-1, \frac{c+1}{2} +N \right)
\end{array}
\]
and
\[
\begin{array}{r}
\displaystyle
\sum_{k=1}^{2N}{{(-1)}^k} \log(k+c)=
\zeta^{\prime}\left(0, \frac{c+1}{2} \right) -
\zeta^{\prime}\left(0, \frac{c+2}{2} \right) +\\
\noalign{\vskip0.5pc} \displaystyle
\zeta^{\prime}\left(0, \frac{c+2}{2} +N \right) -
\zeta^{\prime}\left(0, \frac{c+1}{2} +N \right)
\end{array}
\]
where derivatives $\zeta^{\prime}(-\lambda,z)$ are understood by (\ref{int0}).
Now we use Proposition 1, in particular asymptotic expansions (\ref{asymp}), to get

\begin{equation}
\label{lim}
\begin{array}{r}
\displaystyle
\lim_{N \rightarrow \infty} \sum_{k=1}^{2 N}{(-1)^k} k \,
\Big(\log k - \log(k+2x)\Big) = \\
\noalign{\vskip0.5pc} \displaystyle
2\,x \log \Gamma(x+\frac{1}{2}) - 2\,x  \log \Gamma(x+1) \, - \\
\noalign{\vskip0.5pc} \displaystyle
2\, \zeta^{\prime}(-1, x+\frac{1}{2}) + 2 \, \zeta^{\prime}(-1, x+1)- \\
\noalign{\vskip0.5pc} \displaystyle
\frac{\log 2}{12} - x -3\, \zeta^{\prime}(-1)\\
\end{array}
\end{equation}

\noindent
In view of Proposition 3 with $n=1$, we convert the derivatives of the Hurwitz function to the Barnes function:

\[
\label{toG}
\begin{array}{c}
\displaystyle
\zeta^{\prime}(-1,z) = - \log  G(z+1) + z \log \Gamma (z) + \zeta^{\prime}(-1) \\
\noalign{\vskip1.0pc} \displaystyle

3 \zeta^{\prime}(-1) = -\frac{\log 2}{12} + \frac{\log \pi}{2} +
2\log  G(\frac{1}{2})
\end{array}
\]

\noindent Combining these with (\ref{lim}), after some algebraic
manipulations, leads  to (\ref{prop5}).\qed
\end{pf*}

Here are a few particular cases of (\ref{prop5}):

\begin{equation}
\label{formula1} \lim_{N \rightarrow \infty}
\prod_{k=1}^{2N}{{\big(1+\frac{1}{k}\big)}^{-k (-1)^k}} =
\frac{{A^6}}{e {\sqrt{\pi }}\, 2^{1/6}}
\end{equation}

\begin{equation}
\label{formula2} \lim_{N \rightarrow \infty}
\prod_{k=1}^{2N}{{\big(1+\frac{4}{k}\big)}^{-k (-1)^k}} =
 \frac{3 \,\pi^{2}}{16 \, e^{2}}
\end{equation}

\begin{equation}
\label{formula3} \lim_{N \rightarrow \infty}
\prod_{k=1}^{2N}{{\big(1-\frac{1}{2 k}\big)}^{-k (-1)^k}} =
\frac{A^3 \, e^{-{\bf G}/\pi}  \sqrt{\pi} \, 2^{1/6}}{
\Gamma(\frac{1}{4})}
\end{equation}

\noindent where $A$ is the Glaisher-Kinkelin constant and ${\bf
G}$ is Catalan's constant defined by

\[
\log A =  \frac{1}{12} - \zeta^{\prime}(-1)
\]

\[
{\bf G} = \sum_{k=0}^{\infty}\frac{{{(-1)}^k}}{{{(2 k+1)}^2}}
\]
respectively.

The Melzak product can be further generalized. Here is one of such
formulas

\begin{prop}
The following identity holds for $\Re(x) > - \frac{1}{2}$:

\begin{equation}
\label{prop6}
\begin{array}{l}
\displaystyle
\lim_{N \rightarrow \infty}
\prod _{k=1}^{2N}{{\bigg(1-\frac{4
{x^2}}{{k^2}}\bigg)}^{-{k^2} {{(-1)}^k}}} =
\frac{\cos(\pi  x)}{\pi} \exp\bigg(2 \,{x^2}+\frac{7 \zeta (3)}{2 {{\pi
}^2}}\bigg) \\

\noalign{\vskip0.5pc} \displaystyle

\bigg({G}(1+x){G}(1-x)\bigg)^{4}

{{\Bigg(\frac{{{\Gamma}_3}\big(\frac{3}{2}-x\big)
{{\Gamma}_3}\big(\frac{3}{2}+x\big)}{{{\Gamma}_3}(1-x)
{{\Gamma}_3}(1+x)}\Bigg)}^8}
\end{array}
\end{equation}
\end{prop}
We skip the proof of this proposition, since it could be done in the same manner as in Proposition 5.

\begin{cor}
Identity (\ref{prop6}) can be further simplified to

\begin{equation}
\label{cor}
\begin{array}{c}
\displaystyle
\lim_{N \rightarrow \infty}
\prod _{k=1}^{2N}{{\bigg(1-\frac{4
{x^2}}{{k^2}}\bigg)}^{-{k^2} {{(-1)}^k}}} = \quad \quad \quad \\
\noalign{\vskip0.5pc} \displaystyle

{{\tan(\pi x)}^{-4{x^2}}} \exp\bigg[ 2\,{x^2}(1+\pi i) + \frac{7 \zeta (3)}{2 {{\pi }^2}} \, +\\

\noalign{\vskip0.5pc} \displaystyle

\frac{4 i x}{\pi } \bigg({{\rm Li}_2}(\omega
)-{{\rm Li}_2}(-\omega )\bigg)-\frac{2}{{{\pi
}^2}}\bigg({{\rm Li}_3}(\omega )-
{{\rm Li}_3}(-\omega )\bigg)\bigg]
\end{array}
\end{equation}
where $\omega = \exp(2 \pi i x)$ and ${\rm Li}_{k}(\omega)$ is the polylogarithm, defined by

\[
{\rm Li}_{k}(\omega) = \sum _{k=0}^{\infty} \frac{\omega^{k}}{n^k}, \quad k>1, \> |\omega| \leq 1
\]

\end{cor}

\smallskip
\begin{pf*}{Proof}
The proof follows straightforwardly from (\ref{prop6}) upon
employing the reflexion formulas for the multiple gamma function:

\[
\log \frac{G(1+z)}{G(1-z)} = z \log \left(\frac{\pi}{\sin \pi
z}\right) - \frac{\pi \, i}{2} \,B_{2}(z) + \frac{i}{2 \,\pi}\,
{\rm Li}_{2} (e^{2\pi i z})
\]

\[
\begin{array}{c}
\displaystyle

2\, \log \Big(\Gamma_{3}(1+z) \, \Gamma_{3}(1-z)\Big) + \log
\Big(G(1+z)\, G(1-z) \Big) = \\
\noalign{\vskip 0.5pc} \displaystyle

z^{2} \log \left(\frac{\pi}{\sin \pi z}\right) - \pi \,i \,z
\,B_{2}(z) +
\frac{\pi \, i}{3} \,B_{3}(z)  \, + \\
\noalign{\vskip 0.5pc} \displaystyle

\frac{i \, z}{\pi} \, {\rm Li}_{2} (e^{2\pi i z}) - \frac{1}{2\,
\pi^{2}} \, {\rm Li}_{3} (e^{2\pi i z}) + \frac{\zeta (3)}{2
\pi^{2}}
\end{array}
\]

\noindent where $B_{k}(z)$ are the Bernoulli polynomials.\qed
\end{pf*}

Proposition 6 yields the following particular cases:

\begin{equation}
\label{formula4}
\lim_{N \rightarrow \infty}
\prod_{k=1}^{2N}{{\big(1-\frac{1}{4 {k^2}}\big)}^{-{k^2}
{{(-1)}^k}}} =
\exp \left(\frac{1}{8}-\frac{2 G}{\pi
}+\frac{7 \zeta (3)}{2 {{\pi }^2}}\right)
\end{equation}

\begin{equation}
\label{formula5}
\lim_{N \rightarrow \infty}
\prod_{k=2}^{2N}{{\big(1-\frac{1}{{k^2}}\big)}^{-{k^2}
{{(-1)}^k}}}=\frac{\pi }{4}\exp\left(\frac{1}{2}+\frac{7 \zeta (3)}{{{\pi }^2}}\right)
\end{equation}
Formula (\ref{formula4}) follows from (\ref{cor}) with
$x=\frac{1}{4}$, $\omega=i$ and taking into account

\[
{{\rm Li}_2}(i)-{{\rm Li}_2}(-i)=2 \, i \, {\bf G}
\]

\[
{{\rm Li}_3}(i)-{{\rm Li}_3}(-i)=\frac{i \, {{\pi }^3}}{16}
\]
Similarly, (\ref{formula5}) follows from (\ref{cor})  with
$x=\frac{1}{2}$, $\omega=-1$  and

\[
{{\rm Li}_2}(1)-{{\rm Li}_2}(-1)=\frac{\pi^2}{4}
\]

\[
{{\rm Li}_3}(1)-{{\rm Li}_3}(-1)=\frac{7 \zeta(3)}{4}
\]

\[
\lim_{x \rightarrow \frac{1}{2}} \frac{\tan^{-4 x^2}(\pi x)}{1-4 \,x^2} = \frac{\pi}{4}
\]

\begin{acknowledgements}
This work was supported by grant CCR-0204003 from the National
Science Foundation.
\end{acknowledgements}

\end{document}